\documentclass[12pt]{article}

\usepackage{amsmath}
\usepackage{amssymb}

\newtheorem{theorem}{Theorem}[section]

\newcommand{\PSHG}{{\operatorname{PSHG}}}
\newcommand{\PSH}{{\operatorname{PSH}}}
\newcommand{\PSHS}{{\operatorname{PSHS}}}
\newcommand{\MW}{{\operatorname{MW}}}

\newcommand{\Li}{{\operatorname{L_0}}}
\newcommand{\D}{{\mathbb D}}

\title{Tropical analysis of plurisubharmonic singularities}
\author{Alexander Rashkovskii}
\date{}

\begin{document}

\maketitle

\begin{abstract}
Tropical structures appear naturally in investigation of singularities of plurisubharmonic functions. We show that standard characteristics of the singularities can be viewed as tropicalizations of certain notions from commutative algebra. In turn, such a consideration gives a tool for studying the singularities.
 In addition, we show how the notion of Newton polyhedron and its generalizations come into the picture as a result of the tropicalization.
\end{abstract}

\section{Introduction}
We recall that a semiring $S$ with an addition $\oplus$ and multiplication $\otimes$ is called idempotent if $s\oplus s=s$ for any $s\in S$. When $S$ is a subset of the extended real line,  $a\oplus b=\max\{a,b\}$ (or $\min\{a,b\}$) and $a\otimes b =a+b$, such a semiring is usually called tropical. For basics on idempotent/tropical structures, see, e.g.,
\cite{L} and the bibliography therein.

In this note, we consider certain tropical semirings arising naturally in multidimensional complex analysis. This starts with a simple observation that a basic object of pluripotential theory -- plurisubharmonic functions -- can be viewed as Maslov's dequantization of analytic functions (a basic object of the whole complex analysis). To detect a tropical structure, we need to pass from the world of complex values to the real one. This makes sense in consideration of asymptotic behavior of absolute values $|f(z)|$ of analytic functions $f$ when $z$ approaches either the zero set of $f$ or infinity. Here we will be concerned with the former (local) situation, which invokes investigation of singularities of plurisubharmonic functions and corresponding tropical semirings.

Standard characteristics of singularities of plurisubharmonic functions are thus "tropicalizations" of notions from commutative algebra and can be viewed as functionals on the corresponding tropical semiring. Central role here is played by tropically linear functionals (i.e., additive and multiplicative with respect to the tropical operations and homogeneous with respect to the usual multiplication by positive constants). A problem of description for such functionals is posed. On the other hand, a larger class of the functionals, just tropically additive and positive homogeneous, is described, and a relation between these two classes is established. A few other problems are formulated as well. In addition, we show that the linear functionals can be thought of as "tropicalizations" of valuations on the local ring of germs of analytic functions as well.

Another way of using Maslov's dequantization is to perform it on the arguments of the functions, which moves us from functions on complex manifolds to functions on ${\mathbb R}^n$. This results in a notion of local indicator, introduced from a different point of view in \cite{LeR}. The semiring of plurisubharmonic singularities maps to a tropical semiring of the indicators, and the latter turns out to be isomorphic to an idempotent semiring of complete convex subsets of $\mathbb R_+^n=\{t\in\mathbb R^n:\: t_i> 0,\ i=1,\ldots, n\}$. Going this way, the notion of Newton polyhedron comes naturally into the picture, together with generalizations of famous Kushnirenko's and Bernstein's results on bounds for multiplicities of holomorphic mappings in terms of (mixed) covolumes of the polyhedra.

Most of the results, except for those in Section~\ref{sec:5}, are obtained in \cite{R}--\cite{R7}, so we do not present their proofs here and just put them into the context of tropical mathematics. The proofs of the statements from Section~\ref{sec:5} are sketched.

\section {Plurisubharmonic singularities}
An upper semicontinuous, real-valued function $u$ on a complex manifold $M$
is called {\it plurisubharmonic} ({\it psh}) if for every holomorphic mapping $\lambda$
from the unit disk $\D$ into $M$, the function $u\circ\lambda$ is subharmonic (which means that for every point $\zeta\in\D$, $$(u\circ\lambda)(\zeta)\le \frac{1}{2\pi}\int_0^{2\pi}(u\circ\lambda)(\zeta+re^{i\theta})\,d\theta$$ for all $r<1-|\zeta|$). A basic example is
$u=c\log|f|$ with $c>0$ and a function $f$ analytic on $M$. Moreover, as
follows from Bremermann's theorem \cite{B}, every psh function on a pseudoconvex domain $\omega\subset M$ belongs to the closure (in $L_{loc}^1(\omega)$) of the set of functions $\{\sup_\alpha c_\alpha\log|f_\alpha|\}$. For standard facts on psh functions, see, e.g., \cite{Kl}, \cite{LeG}, \cite{Ro}.

Let ${\cal O}_M$ be the ring of analytic functions $f$ on $M$. The transformation
$f\mapsto\log|f|$ maps it to the cone $\PSH(M)$ of psh functions on $M$, and the ring operations on ${\cal O}_M$ induce a natural tropical structure on $\PSH(M)$ with the
addition $$u\check\oplus v:=\max\{u,v\},$$ which is based on
Maslov's dequantization $$f+g\mapsto \frac1N\log|f^N+g^N|\longrightarrow
\max\{\log|f|,\log|g|\}\quad {\rm as\ }N\to\infty,$$ and
multiplication $$u\otimes v:=u+v$$ (simply by $fg\mapsto \log|f
g|=\log|f|\otimes\log|g|$). Thus $\PSH(M)$ becomes a tropical
semiring, closed under (usual) multiplication by positive constants.
(We use the symbol $\check\oplus$ instead of $\oplus$ in order to
emphasize that it is the max-addition; later on we will need another
idempotent addition, $a\hat\oplus b=\min\{a,b\}$.) The neutral
element (tropical $0$) is $u \equiv-\infty$, and the unit (tropical
$1$) is $u\equiv 0$.

From now on, we restrict ourselves to local considerations, so in the sequel we deal with functions defined near $0\in{\mathbb C}^n$. Let ${\cal O}_0$ denote the ring of germs of analytic functions at $0$, and let $\mathfrak{m}_0=\{f\in {\cal O}_0 :
f(0)=0\}$ be its maximal ideal. The above log-transformation sends
${\cal O}_0$ to the corresponding tropical semiring $\PSHG_0$ of germs of
psh functions. We will say that a psh germ $u$ is {\it
singular} at $0$ if $u$ is not bounded (below) in any neighbourhood of
$0$. For functions $u=\log|f|$ this means $f\in\mathfrak{m}_0$;
asymptotic behaviour of arbitrary psh functions can be much more
complicated (it may even happen that $u(0)>-\infty$). 

A partial order on $\PSHG_0$ is given as follows: $$u\preceq v
\Leftrightarrow u(z)\le v(z)+O(1),\quad z\to 0,$$ which leads to the
equivalence relation $u\sim v$ if $u(z)=v(z)+O(1)$. The equivalence
class ${\rm cl}(u)$ of $u$ is called the {\it plurisubharmonic
singularity} of $u$ (in \cite{Zah}, a closely related object was introduced under the name {"standard singularity"}). The collection of psh singularities
$\PSHS_0=\PSHG_0/\sim$ has the same tropical structure
$\{\check\oplus,\otimes\}$ and the partial order ${\rm cl}(u)\le
{\rm cl}(v)$ if $u\preceq v$. The neutral element here is still
$u\equiv -\infty$, while the unit is represented by any nonsingular
germ.

Psh singularities form a cone with respect to the usual
multiplication by positive numbers. Finally, they are endowed
with the topology where ${\rm cl}(u_j)\to {\rm cl}(u)$ if there
exists a neighbourhood $\omega$ of $0$ and psh functions $v_j\in
{\rm cl}(u_j)$, $v\in {\rm cl}(u)$ such that $v,v_j\in\PSH(\omega)$ and $v_j\to v$
in $L^1(\omega)$. By abusing the notation, we will write
occasionally $u$ for ${\rm cl}(u)$.

\medskip
\section{Characteristics of singularities} {\bf 1.} A fundamental
characteristic of an analytic germ $f\in\mathfrak{m}_0$ is its
multiplicity (vanishing order) $m_f$: if $f=\sum P_j$ is the
expansion of $f$ in homogeneous polynomials, $P_j(tz)=t^jP(z)$, then
$m_f=\min\{j: P_j\not\equiv 0\}$.

The corresponding basic characteristic of singularity of $u\in
\PSHG_0$ is its {\it Lelong number}
$$ \nu(u)=\lim_{t\to -\infty}t^{-1} M(u,t)=
\liminf_{z\to 0}\frac{u(z)}{\log|z|}=dd^cu\wedge
(dd^c\log|z|)^{n-1}(0);$$ here $M(u,t)$ is the mean value of $u$
over the sphere $\{|z|=e^t\}$; $d=\partial + \bar\partial$, $d^c= (
\partial -\bar\partial)/2\pi i$. If $f\in\mathfrak{m}_0$, then
$\nu(\log|f |)=m_f$. This characteristic of singularity gives
important information on the asymptotics:
$u(z)\le\nu(u)\log|z|+O(1)$.

Since $\nu(v)=\nu(u)$ for all $v\in {\rm cl}(u)$, Lelong number can
be considered as a functional on $\PSHS_0$ with values in the
tropical semiring $\mathbb R_{+}(\min,+)$ of non-negative real
numbers with the operations $$x\hat\oplus y=\min\{x,y\} {\quad \rm
and\quad } x\otimes y=x+y.$$ As such, it is
\par\noindent (i) positive homogeneous: $\nu(cu)=c\,\nu(u)$ for all
$c>0$, \par\noindent (ii) additive: $\nu(u\check\oplus v)=
\nu(u)\hat\oplus \nu(v)$, \par\noindent (iii) multiplicative:
$\nu(u\otimes v)= \nu(u)\otimes \nu(v)$, and \par\noindent (iv)
upper semicontinuous: $\nu(u)\ge\limsup \nu(u_j)$ if $u_j\to u$.

\medskip
{\bf 2.} Lelong numbers are independent of the choice of coordinates
(Siu's theorem). Let us now fix a coordinate system centered at
$0$. The {\it directional Lelong number} of $u\in\PSHG_0$ in a
direction $a\in\mathbb R_+^n$ (introduced by C.~Kiselman
\cite{Kis2}) is
\begin{equation}\label{dir} \nu(u,a)=\lim_{t\to -\infty}t^{-1}
M(u,ta)= \liminf_{z\to 0}\frac{u(z)}{\phi_a(z)},\end{equation} where
$M(u,ta)$ is the mean value of $u$ over the distinguished boundary
of the polydisk $\{|z_k|<\exp(ta_k)\}$ and
\begin{equation}\label{phi}
\phi_a(z)=\check\oplus_k \,a_k^{-1}\log|z_k|.\end{equation} Since
its value is constant on $cl(u)$, it is well defined on $\PSHS_0$.
The functional has the same properties (i)--(iv), and the collection
$\{\nu(u,a)\}_{a}$ gives refined information on the singularity $u$;
in particular, $\nu(u)=\nu(u,(1,\ldots,1))$.

For polynomilas or, more generally, analytic functions $f=\sum c_Jz^J\in\mathfrak{m}_0$,
it can be computed as
$$\nu(\log|f|,a)=\hat\oplus\{\langle a,J\rangle:c_J\neq 0\},$$ the expression in the right-hand side being known in number theory
as the {\it index} of $f$ with respect to the weight $a$ \cite{La}.

\medskip
{\bf 3.} A general notion of Lelong number with respect to a
plurisubharmonic weight was introduced by J.-P.~Demailly \cite{D1} (concerning the complex Monge-Amp\`ere operator $(dd^c)^n$, the reader can consult \cite{Kl} and \cite{D}).
Let $\varphi\in \PSHG_0$ be continuous and locally bounded outside
$0$. Then the mixed Monge--Amp\`ere current $dd^c u\wedge
(dd^c\varphi)^{n-1}$ is well defined for any psh function $u$ and is
equivalent to a positive Borel measure. Its mass at $0$,
\begin{equation}\label{eq:ldn}\nu(u,\varphi)= dd^c
u\wedge (dd^c\varphi)^{n-1}(\{0\}),\end{equation} is the {\it generalized}, or
{\it weighted, Lelong number} of $u$ with respect to the weight
$\varphi$. By Demailly's comparison theorem, it is constant on ${\rm
cl}(u)$ and thus defines a functional on $\PSHS_0$. It still has the
above properties (i), (iii), and (iv), however in general is only
subadditive: $\nu(u\check\oplus v, \varphi)\le \nu(u,
\varphi)\hat\oplus \nu(v, \varphi)$.

\medskip
{\bf 4.} One more characteristic, the {\it integrability index}
\begin{equation}\label{ii} \lambda_u=
\inf\{\lambda>0:e^{-u/\lambda}\in L^2_{loc}\},\end{equation} is both
subadditive and submultiplicative, and it is also upper
semicontinuous \cite{DK}. If
$f=(f_1,\ldots,f_m)\in\mathfrak{m}_0^m$, the value
$\lambda_{\log|f|}$ is known as the {\it Arnold multiplicity} of the
ideal ${\mathcal I}$ generated by $f_j$, and
\begin{equation}\label{lc}
lc({\mathcal I})=\lambda_{\log|f|}^{-1}\end{equation} is the {\it
log canonical threshold} of ${\mathcal I}$.

\section{Additive functionals} Another generalization of the notion
of Lelong number was introduced in \cite{R7}. Let $\varphi\in
\PSHG_0$, singular at $0$, be locally bounded and {\it maximal} outside $0$ (that is,
satisfies $(dd^c\varphi)^n=0$ on a punctured neighbourhood of $0$);
the collection of all such germs ({\it maximal weights}) will be
denoted by $\MW_0$. An important example of such a weight is $\varphi=\log|F|$ for an equidimensional holomorphic mapping $F$ with isolated zero at the origin.

For $u\in \PSHG_0$ (or $u\in \PSHS_0$), its {\it
type relative to} $\varphi\in \MW_0$ is defined as
$$\sigma(u,\varphi)=\liminf_{z\to 0}\frac{u(z)}{\varphi(z)}.$$
When $u=\log|f|$ and $\varphi=\log|F|$, the relative type $\sigma(u,\varphi)$
equals the value $\bar\nu_{\mathcal I}(f)$ considered in \cite{LT}, ${\mathcal I}$ being the ideal in ${\cal O}_0$ generated
by the components of the mapping $F$. For the directional
weights $\phi_{a}$ (\ref{phi}),
$$\sigma(u,\phi_a)=\nu(u,a)=a_1\ldots a_n\,\nu(u,\phi_a).$$

Since the function $t\mapsto \sup\,\{u(x):\:\varphi(x)<t\}$ is
convex, one has the relation $$u\preceq
\sigma(u,\varphi)\varphi.$$

Given $\varphi\in \MW_0$, the functional $\sigma(\cdot,\varphi):
\PSHS_0\to [0,+\infty]$ is positive homogeneous, additive,
supermultiplicative, and upper semicontinuous. Actually, relative
types give a general form for all "reasonable" additive functionals
on $\PSHS_0$:

\medskip
\begin{theorem}\label{theo:1} {\rm\cite{R7}} { Let a functional $\sigma:\:
\PSHS_0\to [0,\infty]$ be such that
\begin{enumerate} \item[1)] $\sigma(cu)=c\,\sigma(u)$ for all $c>0$;
\item[2)] $\sigma(\check\oplus u_k)=\hat\oplus \sigma(u_k)$, $k=1,2$;
\item[3)] if $u_j\to u$, then $\limsup\,
\sigma(u_j)\le\sigma(u)$;
\item[4)] $\sigma(\log|z|)=1$;
\item[5)] $\sigma(u)< \infty$ if
$u\not\equiv-\infty$.
\end{enumerate}
Then there exists a weight $\varphi\in \MW_0$ such that
$\sigma(u)=\sigma(u,\varphi)$ for every singularity $u\in \PSHS_0$.
The representation is essentially unique: if two maximal weights
$\varphi$ and $\psi$ represent $\sigma$, then ${\rm
cl}(\varphi)={\rm cl}(\psi)$.}\end{theorem}

In particular, such a functional $\sigma$ is always
supermultiplicative; if $n=1$, it is multiplicative (equal to the
mass of the Riesz measure $\frac1{2\pi}\Delta u$ at $0$).

The function $\varphi\in \MW_0$ from the theorem can be constructed
by the Perron method (as a corresponding Green function): given a
bounded hyperconvex neighbourhood $\Omega$ of $0$, $\varphi$ is the
upper envelope of all negative psh functions $v$ in $\Omega$ such
that $\sigma(v)\ge 1$.

\section{Additivity vs linearity}\label{sec:5} A functional on $\PSHS_0$ is
{\it (tropically) linear} if it is both additive and multiplicative;
we will also assume it to be positive homogeneous and upper
semicontinuous. The collection of all such functionals will be
denoted by $\Li$.

An example of linear functional is $u\mapsto\nu(u\circ \mu,p)$, the
Lelong number of the pullback of $u$ by a holomorphic mapping $\mu$
at $p\in\mu^{-1}(0)$. Another example are the directional Lelong numbers
$\nu(u,a)$ defined by (\ref{dir}).

In a usual vector space, every convex function is an upper envelope
of affine ones. In our situation, any tropically additive functional
is superadditive in the usual sense, while tropically linear ones
are additive. This raises the following question.

\medskip\noindent
{\bf Problem 1.} Is it true that all tropically additive functionals
on $\PSHS_0$ can be represented as lower envelopes of tropically
linear ones?

\medskip

We can answer the question for the functionals generated by
weights that can be approximated by multiples of logarithms of moduli of holomorphic mappings.
First, let us take $\varphi=\log|F|\in \MW_0$, where $F$ is an equidimensional holomorphic mapping with
isolated zero at the origin. By the Hironaka desingularization theorem, there exists a "log resolution" for the mapping $F$, i.e., a proper holomorphic mapping $\mu$ of a manifold $X$ to a neighborhood $U$ of $0$, that is an isomorphism between $X\setminus\mu^{-1}(0)$ and $U\setminus\{0\}$, such that $\mu^{-1}(0)$ is a normal crossing divisor with components $E_1,\ldots,E_N$, and in local coordinates centered at a generic point $p$ of a nonempty intersection $ E_I=\cap_{i\in I} E_{i}$, $I\subset \{1,\ldots,N\}$, $$(F\circ\mu)(x)=h(x)\prod_{i\in I}x_i^{m_i}$$
with $h(0)\neq 0$. Then for any $u\in\PSHS_0$, one has
$$\sigma(u,\varphi)=\min \{\nu_{I,m_I}(u\circ\mu): E_I\neq\emptyset\},$$ where
$$ \nu_{I,m_I}(u\circ\mu)=
\liminf_{x\to 0}\frac{(u\circ\mu)(x)}{\sum_{i\in I} m_i\log|x_i|}$$
at a generic point of $p\in E_I$.
It is then easy to see that
$$\nu_{I,m_I}(u\circ\mu)\ge \min_{i\in I}\nu_{i,m_i}(u\circ\mu)=\min_{i\in I}m_i^{-1}\nu_i(u\circ\mu),$$
where $\nu_i(u\circ\mu)$ is the Lelong number of $u\circ\mu$ at a generic point of $E_i$, which is a linear functional. This gives us the following result (proved for ideals ${\cal I}\subset{\cal O}_0$ in \cite{LT}).

\begin{theorem}\label{theo:anue} For any weight $\log|F|\in\MW_0$ there exist
finitely many functionals $l_j\in \Li$ such that
$\sigma(u,\log|F|)=\min_j\, l_j(u)$ for every $u\in\PSHS_0$; in
other words,
$$\sigma(u,\log|F|)=\hat\oplus_j\:  l_j(u),\quad u\in\PSHS_0.$$
\end{theorem}

Furthermore, let us say that a function $\varphi\in\PSHG_0$ has {\it asymptotically analytic singularity} if for any $\epsilon>0$ there exist positive integers $p$ and $q$, a constant $C>0$, a neighbourhood $U$ of $0$, and a holomorphic mapping $f:U\to{\mathbb C}^p$ such that
\begin{equation}\label{eq:asan}(1+\epsilon)\varphi(z)-C\le q^{-1}\log|f(z)|\le (1-\epsilon)\varphi(z)+C,\quad z\in U.\end{equation}

It can be easily shown that any weight $\varphi\in \MW_0$ with asymptotically analytic singularity has a continuous representative $\psi\in cl(\varphi)\cap\MW_0$ which can be approximated as in (\ref{eq:asan}) with $p=n$ for all $\epsilon>0$. By using Demailly's approximation theorem \cite{D2}, it was shown in \cite{BFaJ1} that
(\ref{eq:asan}) holds if $e^\varphi$ is H\"older continuous or, more generally, if $\varphi$ is a {\it tame} weight, which means that there exists a constant $C_\varphi>0$ such that for any
$t>C_\varphi$ the condition $|f|\exp\{-t\varphi\}\in L^2_{loc}$ for a
function $f\in{\cal O}_0$ implies $\sigma(\log|f|,\varphi)\ge t-C_\varphi$.
(Actually, we are unaware of any example
of maximal weight whose singularity is not asymptotically analytic.)

The following result is a direct consequence of Theorem~\ref{theo:anue}; for tame weights it is essentially proved in \cite{BFaJ1}.

\begin{theorem}  If $\varphi\in\MW_0$ has asymptotically analytic singularity, then
$$\sigma(u,\varphi)=\hat\oplus\,\{l(u):\: l\in \Li,\  l\ge
\sigma(\cdot,\varphi)\},\quad u\in\PSHS_0.$$
\end{theorem}

In view of Theorem~\ref{theo:1}, the following problems are natural.

\medskip\noindent
{\bf Problem 2.} Describe all $\varphi\in\MW_0$ such that the
functional $\sigma(\cdot,\varphi)\in \Li$.

\medskip\noindent
{\bf Problem 3.} What are functionals satisfying all the conditions
of Theorem~\ref{theo:1} except for the last one?

\medskip\noindent
{\bf Problem 4.} Does there exist a functional $\sigma\neq 0$
satisfying conditions 1)--3) and 5) of Theorem~\ref{theo:1}, such
that $\sigma(\log|z|)=0$?

\medskip\noindent
{\bf Problem 5.} What are multiplicative functionals on $\PSHS_0$?

\section{Relative types and valuations} For basics on valuation
theory, we refer to \cite{ZS}.
  Recall that a {\it valuation} on
the analytic ring ${\cal O}_0$  is a nonconstant function $\mu:\:
{\cal O}_0\to [0,+\infty]$ such that
$$\mu(f_1f_2)=\mu(f_1)+\mu(f_2),\quad
\mu(f_1+f_2)\ge\min\{\mu(f_1),\mu(f_2)\}, \quad\mu(1)=0;$$ a
valuation $\mu$ is {\it centered} if $\mu(f)>0$ for all
$f\in\mathfrak{m}_0$.

Every $\varphi\in \MW_0$ generates a functional
$\sigma_\varphi(f)=\sigma(\log|f|,\varphi)$ on ${\cal O}_0$
satisfying
$$\sigma_\varphi(f_1f_2)\ge \sigma_\varphi(f_1)+\sigma_\varphi(f_2), \quad
\sigma_\varphi(f_1+f_2)\ge\min\{\sigma_\varphi(f_1),\sigma_\varphi(f_2)\},
\quad\sigma_\varphi(1)=0.$$ (In \cite{LT} such functions are called {\it
order functions}, and in \cite{R} -- {\it filtrations}.) If the relative type functional
$\sigma(\cdot,\varphi)$ is multiplicative, then $\sigma_\varphi$ is
a valuation, centered if $\sigma(\log|z|,\varphi)>0$.

One can thus consider tropically linear
functionals on $\PSHS_0$ as tropicalizations of certain valuations
on ${\cal O}_0$. For example, the (usual) Lelong number is the
tropicalization of the multiplicity valuation $m_f$. The types
relative to the directional weights $\phi_{a}$ (\ref{phi}) are
multiplicative functionals on $\PSHS_0$, and $\sigma_{\phi_{a}}$ (Kiselman's directional Lelong numbers) are
monomial valuations on ${\cal O}_0$.

It was shown in \cite{FaJ} for $n=2$ and in \cite{BFaJ1} in the
general case that an important class of valuations ({\it quasi-monomial}
valuations) can be realized as
$\sigma_\phi$; all other
centered valuations are limits of increasing sequences of the
quasi-monomial ones.
In addition, the Demailly's weighted Lelong number $\nu(\cdot,\varphi)$ (\ref{eq:ldn}) with a tame weight $\varphi$ is an average
of valuations \cite{BFaJ1}.

\section{Local indicators as Maslov's dequantizations}\label{sec:MD}
Singular psh germs appear as Maslov's dequantization of analytic
functions. As indicated by constructions in tropical algebraic
geometry \cite{V}, it is reasonable to perform a dequantization in
the argument as well. This turns out to be equivalent to
consideration of {\it local indicators}, a notion introduced in
\cite{LeR} by a completely different argument.

For a fixed coordinate system at $0$, let $\nu(u,a)$ be the
directional Lelong number of $u\in \PSHS_0$ in the direction $a\in
\mathbb R_+^n$, see (\ref{dir}). Then the function
$$\psi_u(t)=-\nu(u,-t),\quad
t\in\mathbb R_-^n=-\mathbb R_+^n,$$ 
is convex and increasing in each $t_k$, so
$\psi_u(\log|z_1|,\ldots,|z_n|)$ can be extended (in a unique way)
to a function $\Psi_u(z)$ plurisubharmonic in the unit polydisk
$\D^n$, the {\it local indicator} of $u$ at $0$ \cite{LeR}. Note
that the map $u\mapsto\Psi_u$ keeps the tropical structure:
$$\Psi_{u\check\oplus v}=\Psi_u \check\oplus \Psi_v, \quad
\Psi_{u\otimes v}=\Psi_u \otimes \Psi_v, \quad
\Psi_{c\,u}=c\,\Psi_u.$$

It is easy to see that the indicators have the log-homogeneity
property
\begin{equation}\label{hom}\Psi_u(z_1,\ldots,z_n)=\Psi_u(|z_1|,
\ldots,|z_n|)= c^{-1}\Psi_u(|z_1|^c,\ldots,|z_n|^c) \quad \forall
c>0.\end{equation} In particular, this implies $(dd^c\,\Psi_u)^n=0$
on $\{\Psi_u>-\infty\}$, so if $\Psi_u$ is locally bounded outside
$0$, then $\Psi_u\in \MW_0$,
$$(dd^c\,\Psi_u)^n=N_u\delta_0$$ for some $N_u\ge 0$, and $N_u=0$ if
and only if $\Psi_u\equiv 0$ ($\delta_0$ is Dirac's
$\delta$-function at $0$). 

The indicator can be viewed as the tangent (in the logarithmic
coordinates) for the function $u$ at $0$ in the following sense.

\begin{theorem} {\rm\cite{R}} {The indicator $\Psi_u(z)$ is a unique
$L_{loc}^1$-limit of the functions \begin{equation}\label{um}
T_mu(z)=m^{-1}u(z_1^m,\ldots, z_n^m),\quad
m\to\infty.\end{equation}}
\end{theorem}

In the tropical language, this means that for $f\in \mathfrak{m}_0$,
the sublinear function $\psi_{\log|f|}(t)$ on $\mathbb R_-^n$ is
just a Maslov's dequantization of $f$:
$$\psi_{\log|f|}(t)=\lim_{m\to\infty}m^{-1}\log|f(e^{m(t_1+i\theta_1)},
\ldots, e^{m(t_n+i\theta_n)})|,$$
an interesting moment here being that the arguments become real by themselves.

The indicators are psh characteristics of psh singularities:
\begin{equation}\label{bound} u(z)\le \Psi_u(z)+O(1);\end{equation}
if  $\Psi\in L_{loc}^\infty(\D^n\setminus\{0\})$, then $\Psi_u\in \MW_0$ and the relative type $\sigma(u,\Psi_u)=1$.
When $u$ has isolated singularity at $0$, this implies (by
Demailly's comparison theorem \cite{D1}) the following relation between the
Monge-Amp\`ere measures:
$$(dd^cu)^n\ge (dd^c\,\Psi_u)^n=N_u\delta_0;$$
note that the measures $(dd^c\,T_mu)^n$ of $T_mu$ (\ref{um}) need
not converge to $(dd^c\,\Psi_u)^n$.

More generally, if for an $n$-tuple of psh functions $u_k$ the
current $dd^cu_1\wedge\ldots\wedge dd^cu_n$ is well defined near
$0$, then
$$dd^cu_1\wedge\ldots\wedge dd^cu_n\ge
dd^c\Psi_{u_1}\wedge\ldots\wedge
dd^c\Psi_{u_n}=N_{\{u_k\}}\delta_0.$$

In addition, relation (\ref{bound}) gives an upper bound for the
integrability index $\lambda_u$ (\ref{ii}),
\begin{equation}\label{iii} \lambda_u\ge
\lambda_{\Psi_u};\end{equation} unlike the situation with the
Monge-Amp\`ere measures, one in fact has $\lambda_{T_mu}\to\lambda_{\Psi_u}$
for the functions $T_mu$ defined by (\ref{um}), which follows from a
semicontinuity property for the integrability indices proved in
\cite{DK}.

In the case of a multicircled singularity
$u(z)=u(|z_1|,\ldots,|z_n|)$, one has actually the equalities
$(dd^cu)^n(0)=(dd^c\Psi_u)^n(0)$ (proved in \cite{R4}) and
$\lambda_u= \lambda_{\Psi_u}$, which follows, by the same
semicontinuity property, from the observation that in this case,
$\Psi_u$ is the upper envelope of negative psh functions $v$ in
$\D^n$ such that $v\le u+O(1)$ near $0$.

\section{Indicators and Newton polyhedra}

Since $\Psi_u$ is log-homogeneous, one can compute explicitly its
Monge-Amp\`ere mass and integrability index.

By transition from the psh function $\Psi_u$ to the convex function
\begin{equation} \label{psi}
\psi_u(t)=\Psi_u(e^{t_1},\ldots,e^{t_n}), \quad t\in\mathbb
R_-^n,\end{equation} and from the complex Monge-Amp\`ere operator to
the real one, we get a representation of the Monge-Amp\`ere measures
in terms of euclidian volumes. Let $\langle a,b \rangle$ stand for
the scalar product in ${\mathbb R}^n$. The convex image $\psi_u$ of
the indicator $\Psi_u$ coincides with the support function to the
convex set
$$\Gamma_u=\{b\in\mathbb R_+^n:\: \psi_u(t)\ge\langle b,t \rangle\
\forall t\in\mathbb R_-^n\}=\{b\in\mathbb R_+^n:\:
\nu(u,a)\le\langle a,b \rangle\ \forall a\in\mathbb R_+^n\},$$  that
is, $$\psi_u(t)=\sup\,\{\langle t, a\rangle: \: a\in \Gamma_u\}.$$
These transformations define an isomorphism between the semiring of
the indicators $\Psi$ and the semiring of complete convex subsets
$\Gamma$ of $\mathbb R_+^n$ (the completeness being in the sense
$a\in\Gamma\Rightarrow a+\mathbb R_+^n\in\Gamma$), endowed with the
operations
$$\Gamma_1\dot\oplus\Gamma_2={\rm conv}\,(\Gamma_1\cup\Gamma_2)$$
(the convex hull of the union) and
$$\Gamma_1\dot\otimes\Gamma_2=\Gamma_1+\Gamma_2=\{a+b:\:
a\in\Gamma_1,\ b\in\Gamma_2\}$$ (Minkowski's addition), and
multiplication by positive scalars $c$. We get then
$$\Gamma_{u\check\oplus v}=\Gamma_u \dot\oplus \Gamma_v, \quad
\Gamma_{u\otimes v}=\Gamma_u \dot\otimes \Gamma_v, \quad
\Gamma_{c\,u}=c\,\Gamma_u.$$

Let ${\rm Covol}\,(\Gamma)$ denote the euclidian volume of $\mathbb
R_+^n\setminus\Gamma$.

\begin{theorem}\label{theo:b} {\rm\cite{R}} {The residual
Monge-Amp\`ere mass of $u\in \PSHG_0$ with isolated singularity at
$0$ has the lower bound}
\begin{equation}\label{b1} (dd^cu)^n(0)\ge (dd^c\,\Psi_u)^n(0)= n!\,{\rm
Covol}\,(\Gamma_u).\end{equation}\end{theorem}

Similarly, the mass of the mixed Monge-Amp\`ere current
$dd^c\Psi_{u_1}\wedge\ldots\wedge dd^c\Psi_{u_n}$ (when well
defined) equals $n!\,{\rm Covol}\,(\Gamma_{u_1},
\ldots,\Gamma_{u_n})$, where ${\rm Covol}\,(A_1,\ldots,A_n)$ is an
$n$-linear form on convex subsets of $\mathbb R_+^n$ such that ${\rm
Covol}\,(A,\ldots,A)={\rm Covol}\,(A)$. This gives the relation
$$ dd^c u_1\wedge\ldots\wedge dd^c u_n (0) \ge n!\,{\rm Covol}\,(\Gamma_{u_1},
\ldots,\Gamma_{u_n}),$$ provided the left-hand side is well defined.

\medskip
If $0$ is an isolated zero of a holomorphic mapping $F=(f_1,\ldots,
f_n)$, then its multiplicity equals
$m_F=(dd^c\log|F|)^n(0)=dd^c\log|f_1|\wedge\ldots\wedge
dd^c\log|f_n|(0)$ and the set $\Gamma_{\log|F|}$ is the convex hull
of the union of the {\it Newton polyhedra} $$\Gamma_{\log|f_j|}={\rm
conv} \{J+\mathbb R_+^n: \: D^{(J)}f_j(0)\neq 0\}$$ of $f_j$ at $0$,
$1\le j\le n$. Therefore, Theorem~\ref{theo:b} implies Kushnirenko's
bound \cite{Ku1}
\begin{equation}\label{b10} m_F\ge
n!\,{\rm Covol}\,(\Gamma_{\log|F|}),\end{equation} while the
relation $m_F\ge n!\,{\rm
Covol}\,(\Gamma_{\log|f_1|},\ldots,\Gamma_{\log|f_n|})$ is a
modification of the local variant of D.~Bernstein's theorem
\cite[Theorem~22.10]{AYu}. As shown in \cite{R7}, an equality in
(\ref{b10}) is true if and only if $\log|F|=\Psi_{\log|F|}+O(1)$.

\medskip
The class of all log-homogeneous psh weights $\Psi$ is generated by
the directional weights $\phi_a$ (\ref{phi}), in the sense that
$\Psi(z)=\check\oplus\{\phi_a(z):\phi_a\le \Psi\}$ and the relative
type
\begin{equation}\label{reps}\sigma(u,\Psi)=\hat\oplus\{\nu(u,a):a\in
A_\Psi\},\end{equation} where $A_\Psi=\{a\in\mathbb
R_+^n:\nu(\Psi,a)\ge1\}$. Moreover, the generalized Lelong number
with respect to any log-homogeneous weight can be represented in
terms of the directional numbers, too:

\begin{theorem}\label{theo:genl}
{\rm\cite{R4}} { For each $\varphi\in \MW_0$ there
exists a positive Borel measure $\gamma_\varphi$ on the set
$A_\varphi$ such that}
\begin{equation}\label{b2}
\nu(u,\varphi)\ge
\nu(u,\Psi_\varphi)=\int_{A_\varphi}\nu(u,a)\,d\gamma_\varphi(a),
\quad u\in \PSHS_0. \end{equation}\end{theorem}

Note that representation (\ref{b2}) for $\nu(u,\Psi)$ is a
(tropically) multiplicative analogue of the additive representation
(\ref{reps}) for $\sigma(u,\Psi)$.

The function $\psi_u$ (\ref{psi}) can be considered as the
restriction of the valuative transform of $u$ (action of $u$ on
valuations \cite{FaJ}, \cite{BFaJ1}) to the set of all monomial
valuations. Although relations (\ref{b1}) and (\ref{b2}) are coarser
than the corresponding estimates in terms of the valuative
transforms from \cite{FaJ} and \cite{BFaJ1}, they give bounds that can
be explicitly computed (the measure $\gamma_\varphi$ is defined
constructively \cite{R4}). This reflects one of the benefits of
using tropical mathematics.

\medskip
Finally, a direct computation involving the function $\psi_u$
(\ref{psi}) shows that the integrability index (\ref{ii}) for the
indicator $\Psi_u$ can be computed as
$\lambda_{\Psi_u}=\sup\{\psi_u(t)/\sum t_j\}$. By (\ref{iii}), it
gives the bound in terms of the directional numbers $\nu(u,a)$:
$$\lambda_{u}\ge\lambda_{\Psi_u}=\sup\{\nu(u,a):\:\sum_k a_k=1\},$$
with an equality in the case of multicircled singularity $u$ (see
the remark in the end of Section~\ref{sec:MD}). This recovers
\cite[Thm. 5.8]{Kis3}, which in turn implies a formula for the log
canonical threshold (\ref{lc}) for monomial ideals proved
independently in \cite{H}.

\section{Related topics} The results on local indicators have
global counterparts for psh functions of logarithmic growth in
$\mathbb C^n$ (i.e., $\limsup_{|z|\to\infty}u(z)/\log|z|<\infty$, a
basic example being $u=\log|P|$ for a polynomial mapping $P$), see
\cite{R1} and \cite{R2}; they are also connected with the notion of
{\it amoebas} of holomorphic functions \cite{PR}. Similar notions
concerning Maslov's dequantization in $\mathbb C^n$ and generalized
Newton polytops were also considered in \cite{LS}.

\vskip1cm

Tek/Nat, University of Stavanger, 4036 Stavanger, Norway

\vskip0.1cm

{\sc E-mail}: alexander.rashkovskii@uis.no

\end{document}